\documentclass[oneside, 12pt]{amsart}
\usepackage{amscd, amssymb, amsmath, mathrsfs,mathabx}
\usepackage[english]{babel}
\usepackage{booktabs}
\usepackage{tikz-cd}
\usepackage{enumitem}
\usepackage{xypic}
\usepackage[pdftex, colorlinks=true,  citecolor=blue, linkcolor=blue, linktocpage=true]{hyperref} 

\DeclareUrlCommand\arXiv{\urlstyle{same}}

\renewcommand{\labelenumi}{(\roman{enumi})}

\newcommand\mylabel[1]{\label{#1}\marginpar{\vspace{-1ex}\medskip\medskip\footnotesize \tt #1}}
\renewcommand\mylabel[1]{\label{#1}}
\newcommand{\mydate}{
\number\day\space
\ifcase\month \or January\or February\or March\or April\or May\or June\or July\or August\or September\or October\or November\or December\fi 
\space\number\year}

\newtheorem{theorem}{Theorem}[section]
\newtheorem*{maintheorem}{Theorem}
\newtheorem{lemma}[theorem]{Lemma}
\newtheorem{proposition}[theorem]{Proposition}
\newtheorem{corollary}[theorem]{Corollary}

\theoremstyle{definition}
\newtheorem{definition}[theorem]{Definition}
\newtheorem{example}[theorem]{Example}

\newtheorem{remark}[theorem]{Remark}

\newtheorem{para}[theorem]{}
\newtheorem*{acknowledgement}{Acknowledgement}

\theoremstyle{remark}
\newtheorem{claim}[theorem]{Claim}

\newcommand{\eq}[2]{\begin{equation}\label{#1}#2 \end{equation}}

\newcommand{\ZZ}{\mathbb{Z}}

\newcommand{\CC}{\mathbb{C}}

\newcommand{\PP}{\mathbb{P}}
\renewcommand{\AA}{\mathbb{A}}

\newcommand{\shL}{\mathscr{L}}

\newcommand{\catC}{\mathcal{C}}

\newcommand{\alg}{\text{\rm alg}}

\newcommand{\Aut}{\operatorname{Aut}}

\newcommand{\depth}{\operatorname{depth}}

\newcommand{\et}{{\text{\rm et}}}

\newcommand{\Frac}{\operatorname{Frac}}

\newcommand{\Gal}{\operatorname{Gal}}

\newcommand{\Hom}{\operatorname{Hom}}

\newcommand{\lra}{\longrightarrow}
\newcommand{\xr}[1] {\xrightarrow{#1}}

\newcommand{\inj}{\hookrightarrow}
\newcommand{\surj}{\twoheadrightarrow}

\newcommand{\maxid}{\mathfrak{m}}

\renewcommand{\O}{\mathscr{O}}

\newcommand{\op}{\text{\rm op}}

\newcommand{\pr}{\operatorname{pr}}

\newcommand{\ra}{\rightarrow}

\newcommand{\red}{{\operatorname{red}}}

\newcommand{\sep}{{\operatorname{sep}}}

\newcommand{\Sch}{\text{\rm Sch}}

\newcommand{\Spec}{\operatorname{Spec}}

\DeclareMathOperator{\FinSet}{FinSet}
\DeclareMathOperator{\FinEt}{FinEt}

\numberwithin{equation}{section}

\begin{document}

\title[]
      {Loops on schemes and \\the algebraic fundamental group}

\author[Kay R\"ulling]{Kay R\"ulling}
\address{Bergische Universität Wuppertal, Faculty of Mathematics and Natural Sciences, 42097 Wuppertal, Germany}
\curraddr{}
\email{ruelling@uni-wuppertal.de}

\author[Stefan Schr\"oer]{Stefan Schr\"oer}
\address{Heinrich Heine University Düsseldorf, Faculty of Mathematics and Natural Sciences, Mathematical Institute, 40204 D\"usseldorf, Germany}
\curraddr{}
\email{schroeer@math.uni-duesseldorf.de}


\dedicatory{Revised version 15.\ February 2024}

\begin{abstract}
 In this note we give a re-interpretation of the algebraic fundamental group for proper   schemes
 that is rather close to the original definition of the fundamental group for topological spaces.
 The idea is to replace the standard interval from topology by what we call interval schemes. This leads to 
 an algebraic version of continuous loops, and the 
 homotopy relation is defined in terms of the monodromy action.
 Our main results hinge on  Macaulayfication for proper schemes  and   Lefschetz type results.
\end{abstract}

\maketitle
\tableofcontents

\section*{Introduction}
\mylabel{introduction}
The fundamental group $\pi_1(X,x_0)$ of a  topological space  $X$ with respect to a base point $x_0$ is an invariant of great significance, 
even more so as its definition is elementary  and intuitive:
the elements are loops up to homotopy, where a loop is a continuous morphism $I\to X$ from the standard interval $I=[0,1]$, 
such that the end points  are mapped to
the base point $x_0$. Roughly speaking, two  loops are homotopic if one can be deformed to the other, respecting the base point.
For connected and locally simply-connected spaces $X$, one may interpret the fundamental group 
also as the group of deck transformations of the universal covering $\tilde{X}\ra X$. 

In the realm of algebraic topology, where one works with schemes rather than topological spaces, the first construction above makes little sense.
However,  Grothendieck  \cite{SGA 1} realized that the second description has an analog in algebraic geometry.
He introduced the notion of a {\em Galois category}  $\catC$, the objects of which should be considered as abstract finite coverings,
which admits a fiber functor $\Phi:\catC\to (\FinSet)$ to the category of finite sets satisfying certain properties.
These properties ensure that   the automorphism group of $\Phi$ is equal to the opposite group of
the automorphisms of an abstract pro-finite universal covering. The algebraic fundamental group $\pi_1^\alg(X,x_0)$ of a connected scheme $X$ with 
respect to a geometric point
$x_0$ is then defined by applying this general construction to the category $(\FinEt/X)$ of finite \'etale coverings of $X$,
with fiber functor given by base-changing along $x_0$. 

If $X$ is of finite type over the complex numbers,
the group $\pi_1^\alg(X,x_0)$ equals the pro-finite completion of the classical fundamental group of 
$X(\CC)$, endowed with the classical topology. If $X=\Spec(K)$ is the spectrum of a field  and $x_0$ is given
by some separable closure, $\pi_1^\alg(X,x_0)$ gives back the corresponding Galois group.
 
 In this  note we observe that the original construction of the fundamental group using loops
 has a meaningful analogy for schemes, once the notions of intervals and loops are
 interpreted in an algebraic manner. More precisely, the crucial properties of the interval $I=[0,1]$ in the construction above 
 are: $I$ is connected, quasi-compact, one dimensional, with no non-trivial coverings,  and is endowed with
 two distinguished points. Translating these properties to algebraic geometry we define in Section  \ref{Intervall schemes}   {\em interval schemes}
as  reduced, connected, affine,    one-dimensional schemes $I$,  which have no non-trivial finite \'etale coverings
and contain two  distinguished closed points with separably closed residue fields. The latter are called \emph{end points}.
An {\em algebraic loop} on a scheme $X$ based at a geometric point $x_0$ is  a morphism of schemes $I\to X$ mapping the end points to the base point.

Algebraic loops define  monodromy transformations. We call   two algebraic loops   \emph{homotopic}
if the resulting  monodromies agree. The {\em algebraic loop group} $\pi_0\Omega^\alg(X,x_0)$ 
is defined as the set of homotopy classes of algebraic loops; the group structure is induced by concatenating
algebraic loops, see Section \ref{Loop space}.

Interval schemes are very often non-noetherian. One example of an interval scheme is the universal Galois covering (introduced by Grothendieck as a pro-object)
of a noetherian, connected, affine, reduced, and one dimensional scheme. 
Such universal Galois coverings were systematically studied in \cite{Vakil Wickelgren 2011}, where Vakil and Wickelgren define  the fundamental group scheme
using  universal coverings, which are certain pro-finite \'etale maps. The notion of interval scheme introduced above is also inspired by their work.
But there are many other  examples of interval schemes, which are more direct to obtain.
For example, if $R$ is  an integral noetherian one-dimensional ring and  $A$ is its integral closure
in the separable closure of $\Frac(R)$, then the choice of two  closed geometric points in $\Spec(A)$
turns it into an interval scheme.

By construction, the monodromy induces an injective   homomorphism
\[\tag{1}
\label{monodromy transformation}
\pi_0\Omega^\alg(X,x_0)^{\op}\lra \pi_1^\alg(X,x_0),
\]
of groups, where ``$\op$'' refers to the opposite group structure.
The main result of this note is the following, see Theorem \ref{thm:loop-et}:

\begin{maintheorem}
Let $X$ be a connected scheme that is separated and of finite type over a ground field $k$,
endowed with a geometric point  $x_0:\Spec (k^\sep)\to X$.
Then the injection \eqref{monodromy transformation} has dense image.
It is actually bijective, provided that $X$ is proper. 
\end{maintheorem}

The main step in the proof of the above theorem for proper $X$  is a Lefschetz type result
saying that for a proper and connected $k$-scheme $X$ we find a closed connected curve $C \subset X$ such that the
algebraic fundamental group of $C$ surjects to the one of $X$, see Proposition \ref{C over a field}. This is well-known 
in the case where $X$ is Cohen--Macaulay and projective over a field, see \cite{SGA 2}, Expos\'e XII. 
We reduce the general situation to this using a van-Kampen-like argument and Macaulayfication,
which was in a special case constructed  by Faltings \cite{Faltings 1978} and in full generality by Kawasaki  \cite{Kawasaki 2000}.
For further results on Macaulayfication, see the recent work of \v{C}esnavi\v{c}ius \cite{Cesnavicius 2021}.
The proof of Theorem \ref{thm:loop-et} is given in Section \ref{sec:proof}.
We do not expect the map $\eqref{monodromy transformation}$ to be an isomorphism for affine schemes of  finite type over a field  in general.

\begin{acknowledgement}
We thank Laurent Moret-Bailly  and the referees for valuable comments, which helped to improve the paper,
and Raju Krishnamoorthy for bringing the reference \cite{Esnault 2017} to our attention.
Furthermore, we thank Oliver Br\"aunling for pointing out that our original title \emph{Algebraic loop groups}
was ambiguous, since the term already refers to certain ind-affine groups attached to linear groups,
compare for example \cite{Faltings  2003}.
This research was conducted in the framework of the   research training group
\emph{GRK 2240: Algebro-Geometric Methods in Algebra, Arithmetic and Topology}, which is funded by the Deutsche Forschungsgemeinschaft.
\end{acknowledgement}

\section{Monodromy}
\mylabel{Monodromy}

Let $Y$ be a scheme, and write $(\FinEt/Y)$ for the category of $Y$-schemes $X$ whose structure
morphism $f:X\ra Y$ is finite and \'etale.
Note that such an $f$ is proper, affine, flat, of finite presentation, and 
for each point $b\in Y$ the fiber $f^{-1}(b)$ is the spectrum of some \'etale algebra
over the residue field $k=\kappa(b)$. 
For the following result, see for example 
\cite{Lenstra 1985}, Theorem 5.10 and Exercise 5.21.

\begin{proposition}\label{prop:split}
Suppose $Y$ is connected, and $f:X\ra Y$ finite and \'etale. Then there is a surjective finite \'etale morphism $Y'\ra Y$
such that $X'=X\times_Y Y'$ is isomorphic over $Y'$ to the disjoint union $\coprod_{i=1}^r Y'$\ 
for some integer $r\geq 0$.
\end{proposition}

This has an immediate consequence:

\begin{corollary}
Suppose $Y$ is connected, and $f:X\ra Y$ finite and \'etale. Then $X$ has only finitely many connected components $U\subset X$,
each of which is open-and-closed. Moreover, the induced morphism $U\ra Y$ is finite and \'etale.
\end{corollary}

\proof
Take $Y'\to Y$ as in Proposition \ref{prop:split}. Since the projection $X'\to X$ is surjective, the connected components
of $X$ are images of the connected components of $X'$, hence there are only finitely many. This implies that
the connected components $U$ of $X$ are open and closed and hence the composition $U\inj X\to Y$ is \'etale and finite.
\qed

\medskip
The proposition tells us that the $Y$-scheme $X$ is a \emph{twisted form} of the disjoint union $\coprod_{i=1}^r Y$,
with respect to the \'etale topology. It thus corresponds to a class 
in the non-abelian cohomology set $H^1(Y,S_r)$, with coefficients in the symmetric group $S_r$
on $r\geq 0$ letters (\cite{Giraud 1971}, Chapter III, Section 2.3). 
To summarize:
\begin{proposition}\label{prop:split2}
If $Y$ is connected, the following are equivalent:
\begin{enumerate}
\item Every finite \'etale $Y$-scheme is isomorphic to some $\coprod_{i=1}^r Y$, $r\geq 0$.
\item We have $H^1(Y,S_r)=\{*\}$ for all integers $r\geq 0$.
\item Each finite \'etale morphism $X\to Y$ from a non-empty connected scheme X is an isomorphism.
\end{enumerate}
\end{proposition}

\begin{definition}\label{defn:simplyconnected}
We say a connected scheme $Y$ is {\em simply connected}, if it satisfies the equivalent conditions of Proposition \ref{prop:split2}.
\end{definition}

\medskip
Let $X$ be a $Y$-scheme, with structure morphism $f:X\ra Y$, and   $a:A\ra Y$ be some other morphism. To simplify notation, we  write
$$
X(A)=\Hom_Y(A,X)=\{a':A\lra X\mid f\circ a'=a\}
$$
of liftings of $a:A\ra Y$ with respect to   $f:X\ra Y$.

\begin{proposition}\label{prop:monodromy}
Suppose that $Y$ is connected, with   $H^1(Y,S_r)=\{*\}$ for all $r\geq 0$, and   $f:X\ra Y$ is finite and \'etale.
Let $a:A\ra Y$ and $b:B\ra Y$ be morphisms with  connected and non-empty domains. Then the sets $X(A)$ and $X(B)$ are finite,
and for each $a'\in X(A)$  there is a unique $b'\in X(B)$ such that the images 
$a'(A)$ and $b'(B) $ lie in the same connected component of $X$.
\end{proposition}
 
\proof
By Proposition \ref{prop:split2}, we may assume $X=\coprod_{i=1}^r Y$, for some $r\ge 0$.
The set $X(A)$ is in bijection  with the set of sections of $\coprod_{i=1}^r A= X\times_Y A\to A$, whence is finite,
and we see that every $a'\in X(A)$ corresponds uniquely to one of the maps
\[A \lhook\joinrel\lra \coprod_i A\xr{\coprod a} \coprod_i Y=X\]
given by including $A$ into one of the $r$ summands. This implies the statement.
\qed

\begin{para}\label{para:monodromy}
In the situation of Proposition \ref{prop:monodromy} we obtain a mapping
$$
\mu_X:X(A)\lra X(B),\quad a'\longmapsto b',
$$
which is called the   \emph{monodromy}, and will play a crucial role throughout.  We   regard it as a natural transformation
between $X\mapsto X(A)$ and $X\mapsto X(B)$, viewed as   functors $(\FinEt/Y)\ra (\FinSet)$.
By Proposition \ref{prop:monodromy}, the monodromy $\mu_X$ is a natural isomorphism, and is given as the composition
of the following natural bijections
\[\mu_X: X(A)\lra\pi_0(X\times_Y A)\lra\pi_0(X)\lra \pi_0(X\times_Y B)\lra X(B),\]
where $\pi_0(X)$ denotes the set of connected components of $X$.
\end{para}

\section{Galois categories}
\mylabel{Galois categories}

In this section we recall the notion of Galois categories, which were introduced
by Grothendieck to unify Galois theory from algebra and covering space theory from topology (\cite{SGA 1}, Expos\'e V).

\begin{para}
Recall that a category $\catC$ is  called a \emph{Galois category}
if there exists a functor $\Phi:\catC\ra (\FinSet)$ such that the following six axioms hold:
\renewcommand{\labelenumi}{(G\arabic{enumi})}
\begin{enumerate}
\item 
Fiber products and final objects exist  in $\catC$.
\item
Finite sums and quotients by finite group actions exist in $\catC$.
\item
Every morphism $X'\ra X$ in $\catC$ factors into a strict epimorphism $X'\ra U$ and the inclusion
of a direct summand $U\subset X$.
\item
The functor $\Phi$ commutes with fiber products and final objects.
\item
It also commutes with finite direct sums and forming quotients by finite group actions, and  transforms strict epimorphisms into surjections.
\item
If for a morphism $u:X'\ra X$ in $\catC$ the resulting  map $\Phi(u)$ is bijective, then $u$ is an isomorphism.
\end{enumerate}
\renewcommand{\labelenumi}{(\roman{enumi})}

One calls $\Phi $ a \emph{fundamental functor} or \emph{fiber functor} for the Galois category $\catC$, and denotes by 
$\pi=\Aut(\Phi)$ the group of natural isomorphisms of the fundamental functor to itself.
In turn, we have an inclusion 
$$
\pi\subset\prod_{X\in\catC} S_{\Phi(X)}
$$
inside a product of   symmetric groups $\Aut(\Phi(X))=S_{\Phi(X)}$.
These groups are finite. We endow them  with the discrete topology, and the product with the product topology.
The latter becomes a topological group that is compact and totally disconnected. 
Such topological groups are also called \emph{pro-finite groups}.
One easily checks  that the subgroup $\pi$ is closed, and thus inherits the structure of a pro-finite group.
Every fiber functor on $\catC$ is (non-canonically) isomorphic to $\Phi$ and hence,
up to a uncanonical isomorphism, the pro-finite group $\pi$ depends only on the Galois category $\catC$, and not
on the choice of fiber functor $\Phi$.

Now write $(\text{$\pi$-FinSet})$ for the category of finite  sets $F$ endowed with a $\pi$-action from the left, 
where the kernel of the canonical homomorphism
$\pi\ra S_F$ is closed. In other words, the action $\pi\times F\ra F$ is continuous, 
when the finite set $F$ is endowed with the discrete topology.
With respect to the natural $\pi$-action on the $\Phi(X)$, $X\in\catC$, the fundamental functor becomes a functor
$$
\Phi:\catC\lra (\text{$\pi$-FinSet}),
$$
and Grothendieck deduced from the axioms (G1)--(G6) that the above is an equivalence of categories.
Conversely, if $G$ is a pro-finite group, the category $(\text{$G$-FinSet})$ is a Galois category:
The functor $\Phi$ that forgets the $G$-action is a fundamental functor, and the resulting $\pi=\Aut(\Phi)$
becomes identified with  $G$. 

One should see $\pi=\Aut(\Phi)$ as a common generalization of the  opposite Galois group $\Gal(F^\sep/F)^\op$ 
for fields $F$, and the pro-finite completion $\widehat{\pi}_1(Y,y_0)$ of the fundamental group, 
say for connected and locally simply-connected topological spaces $Y$.
\end{para}
 
\begin{para}\label{para:functorGalois}
Let  $H:\catC\to\catC'$ be an exact covariant functor between Galois categories. 
By \cite[Expos\'e V, Proposition 6.1]{SGA 1} the exactness is equivalent to the statement, that
the composition  $\Phi'\circ H$ is a fiber functor for $\catC$, whenever $\Phi'$ is a fiber functor for $\catC'$. 
An  exact functor $H$ induces a morphism $h:\pi'\to\pi$ between the corresponding fundamental groups with reversed  direction, 
which is  well-defined up to conjugation. 
The following statements are equivalent (see \cite[Expos\'e V, Proposition 6.9]{SGA 1})
\begin{enumerate}[label=(\arabic*)]
\item $H: \catC\to \catC'$ is fully faithful.
\item $h:\pi'\to \pi$ is surjective.
\item For each connected object $X\in \catC$, the object $H(X)$ is connected. 
\end{enumerate} 
\end{para}

\begin{para}
Let $Y$ be a connected scheme. Then $\catC=(\FinEt/Y)$ becomes a Galois category:
For each morphism $y_0:\Spec(K)\ra Y$, where $K$ is a separably closed field, we get a fiber functor
$$
\Phi_{y_0}:(\FinEt/Y)\lra (\FinSet),\quad X\longmapsto X(K)
$$
given by the set of morphisms $b':\Spec(K)\ra X$ lifting the given $b:\Spec(K)\ra Y$.
The resulting group
$$
\pi_1^\alg(Y,y_0) = \pi= \Aut(\Phi_{y_0})
$$
is called the \emph{algebraic fundamental group} of the connected scheme $Y$ with respect to  $y_0$.
As in  topology, the latter is called   \emph{base point}.
\end{para}

\section{Interval schemes}
\mylabel{Intervall schemes}

In algebraic topology, the \emph{standard interval} $I=[0,1]$ plays a central role.
From our perspective, the following are the crucial properties:
\begin{enumerate}
\item The topological space $I$ is  connected and quasi-compact.
\item There are two distinguished points $0,1\in I$.
\item The universal covering $\tilde{I}\ra I$ is a homeomorphism.
\item The interval $I$ is one-dimensional.
\end{enumerate}

In this section we introduce a class of schemes  with analogous properties.
Fix a separably closed field $K$.

\begin{definition}
\mylabel{interval scheme}
An \emph{interval scheme with $K$-valued endpoints} is a triple $(I,a_0,a_1)$, where $I$   
is a reduced, connected, simply connected, affine, and  one-dimensional scheme and $a_i:\Spec(K)\ra X$  are  two closed embeddings.
\end{definition}

The point  $a_0$ is called the \emph{left endpoint}, whereas  $a_1$ is the \emph{right endpoint}.
By abuse of notation, we usually write $I$ for the interval scheme $(I,a_0,a_1)$. 
The simplest example of an interval scheme is $I=\AA^1_k$ where $k$ is an algebraically closed field of characteristic zero,
and the end points are the rational points $a_0=0$ and $a_1=1$.
However, if $k$ is not algebraically closed or of positive characteristic, this will not be an interval scheme.  
In fact,  interval schemes are very often  non-noetherian.
The following gives the most basic class of interval schemes:

\begin{proposition}
\mylabel{separably closed}
Let $A$ be a one-dimensional integral ring that is normal and whose field of fractions $F=\Frac(A)$ is separably closed, and suppose there
are two surjections $\varphi_i:A\ra K$. Then $I=\Spec(A)$ becomes an interval scheme, where
the endpoints $a_i $ correspond to the homomorphisms $\varphi_i$. 
\end{proposition}

\proof
Let $X\to I$ be a connected finite \'etale covering.  It follows that $X$ is an affine connected normal scheme, and thus  is integral, see \cite{EGAIV4}, Proposition (17.5.7). 
In particular the  function field $L$ of $X$ is a finite separable field extension  $L/\Frac(A)$.
Since $\Frac(A)$ is separably closed we find $L\cong \Frac(A)$. Hence  $X\to I$ is an isomorphism. 
Thus $I$ has no non-trivial finite \'etale covering and therefore defines an interval scheme.
\qed

\begin{example}
Rings as in Proposition \ref{separably closed} easily occur as follows: Suppose that $R$ is a one-dimensional noetherian ring,
endowed with two integral homomorphisms $\psi_i:R\ra K$. The latter means that each $\lambda\in K$
is the root of a monic polynomial with coefficients from  $R$.
Choose a separable closure $F^\sep$ for the field of fractions $F=\Frac(R)$,
and write $A=R^\sep$ for the integral closure of $R\subset F^\sep$.
By construction $A$ is integral and normal, with field of fractions $F^\sep$, and the ring extension $R\subset A$ is integral.
According to the Going-Up Theorem, the map $\Spec(A)\ra \Spec(R)$ is surjective,
thus the $\psi_i:R\ra K$ extend to some homomorphisms $\varphi_i:R^\sep\ra K$. 
Proposition \ref{separably closed} yields that  the scheme $I=\Spec(R^\sep)$ in an interval scheme, where the  endpoints $a_i$
correspond to the homomorphisms $\varphi_i$.
\end{example}

Recall that a local ring $R$ is called \emph{strictly henselian} if each factorization $P\equiv P_1P_2$  into coprime 
polynomials over the residue field $k=R/\maxid_R$  of a monic polynomial $P\in R[T]$ 
is induced by a factorization over $R$, and moreover the   residue field   is 
separably closed. See \cite{EGAIV4}, Proposition (18.8.1) for the next example of an interval scheme, for which the image points of the two end points agree.
\begin{proposition}
\mylabel{strictly henselian}
Let $A$ be a   one-dimensional local ring that is strictly henselian, and whose residue field  
is isomorphic to $K$, and let $\varphi_i:A/\maxid_A\ra K$ be two isomorphisms.
Then $I=\Spec(A)$ becomes an interval scheme, where the endpoints $a_i $ correspond to the homomorphisms $\varphi_i$.
\end{proposition}

%

\begin{para}\label{para:join}
Let $I,J$ be two interval schemes, with $K$-valued endpoints $a_0,a_1$  and $b_0,b_1$, respectively.
We write $I* J$ for the concatenation of $I$ and $J$ with respect to the right endpoint on $I$ and the left endpoint on $J$.
In other words, we have a cocartesian square in the category of schemes
\eq{para:join1}{
\begin{CD}
\Spec(K)    @>b_0>>    J\\
@Va_1VV                                @VVV\\
I        @>>>            I*J.
\end{CD}
}
Note that the cocartesian square above exits in the category of schemes 
by \cite{Ferrand 2003}, Th\'eor\`eme 5.4, and is in fact also cartesian.
The scheme $I*J$ comes with closed embeddings of $I$ and $J$, and we take $a_0\in I\subset I\star J$ as new left endpoint,
and $b_1\in J\subset I* J$ as new right endpoint.
\end{para}

\begin{lemma}\label{lem:join-interval}
In the above situation, the concatenation $I*J$ is an interval scheme, with endpoints $a_0$ and $b_1$.
\end{lemma}

\proof
Write $I=\Spec (A)$ and $J=\Spec (B)$. Then $I* J=\Spec (A\times_ {K} B)$  is affine and reduced. 
Here the fiber product is formed with respect to the homomorphisms $A\rightarrow K\leftarrow B$ corresponding
to the morphisms $a_1$ and $b_0$. 
By construction we have 
$I*J= I\cup J$ and $I\cap J=\Spec(K)$, where the latter can be viewed as the image of both $a_1$ and $b_0$. 
Hence $I*J$ is also one-dimensional  and connected.
Let $X\to I*J$ be a  finite \'etale covering.
Set  $X(Z)=\Hom_{I*J}(Z,X)$, where $Z$ is  an $I*J$-scheme.
By \cite{EGAIV4}, Corollaire (17.9.4), the set $X(I*J)$ is in bijection with the connected components of $X$.
By construction of $I*J$ the set $X(I*J)$ is the pushout of  the maps $X(I)\leftarrow X(I\cap J)\rightarrow X(J)$ induced by
the morphisms $a_1$ and $b_0$. Since $I$, $J$, and $\Spec (K)$ are simply connected the monodromy 
(see \ref{para:monodromy}) induces bijections
\[X(I)\cong X(a_1)\cong X(I\cap J)\cong  X(b_0)\cong X(J).\]
Hence $|X(I*J)|=|X(I\cap J)|=\deg(X_{I\cap J}\to I\cap J)=\deg(X\to I*J)$, 
where $X_{I\cap J}\to I\cap J$ denotes the base change of $X\to I*J$ along $I\cap J\inj I*J$.
Thus any finite \'etale covering of $I*J$ is trivial and hence  $I*J$ is an interval scheme.
\qed

\section{The algebraic loop group}
\mylabel{Loop space}

Fix some separably closed field $K$ and let $Y$ be a scheme, endowed with 
two $K$-valued points  $y_i:\Spec(K)\ra Y$. We may regard this
as an object in the category $(K^2/\Sch)$ of  schemes endowed with two $K$-valued points.

\begin{definition}
\mylabel{algebraic path}
An \emph{algebraic path} in $Y$ starting at $y_0$ and ending in $y_1$
is an interval scheme $(I,a_0,a_1)$ with $K$-valued endpoints, together with a morphism
$$
w:(I,a_0,a_1)\lra (Y,y_0,y_1)
$$
 of  schemes endowed with two $K$-valued points. 
An algebraic path $w$ is called {\em algebraic loop} if $y_0=y_1$.
\end{definition}

By abuse of notation, we often write $w:I\ra Y$ for the algebraic path $w:(I,a_0, a_1)\to (Y, b_0, b_1)$.
For each finite \'etale map $X\ra Y$, the base change induces a finite \'etale map $X\times_Y I\ra I$,
which takes the form $\coprod_{i=1}^r I$, for some $r\geq 0$, and  we have an identification
$X(y_i)= (X\times_Y I)(a_i)$ of fiber sets. In turn, the  monodromy gives a transformation
$$
\mu_w:X(y_0)=(X\times_Y I)(a_0)\lra (X\times_YI)(a_1)=  X(y_1)
$$
that is bijective, and natural in the objects $X\in(\FinEt/Y)$. In other words,
the monodromy $\mu_w$ attached to the path $w:I\ra Y$ from $y_0$ to $y_1$ is a  bijective natural transformation
between   fiber functors 
$$
\Phi_{y_0}, \Phi_{y_1}:(\FinEt/Y)\lra(\FinSet).
$$
We now use this monodromy to give an  algebraic version of homotopy:

\begin{definition}
\mylabel{homotopic}
We say that two algebraic paths $w:I\ra Y$ and $v:J\ra Y$ from $y_0$ to $y_1$
are \emph{homotopic} if $\mu_w=\mu_v$ as  natural transformations 
from $\Phi_{y_0}$ to $\Phi_{y_1}$.
\end{definition}

\begin{para}\label{para:loop-group}
We denote the class of algebraic loops in $Y$ based at the geometric point $y_0$ by
$$
\Omega^\alg(Y,y_0)=\left\{w:(I,a_0,a_1)\lra (Y,y_0,y_0) \mid (I,a_0,a_1) \text{ interval scheme }\right\}.
$$
Let  $w:(I, a_0,a_1)\ra (Y, y_0,y_0)$ and $v:(J, b_0,b_1)\ra (Y,y_0,y_0)$ be two loops based at $y_0$.
It follows from the pushout diagram \eqref{para:join1} and Lemma \ref{lem:join-interval}
that we can concatenate these loops to get a new loop
\[w*v: (I*J, a_0, b_1)\lra (Y,y_0,y_0).\]
The monodromy transformation corresponding to $w*v$ can be factored, for $X\in (\FinEt/Y)$ as 
\[\xymatrix{
X(y_0)\ar[rr]^{\mu_{w*v}} & &  X(y_0) \\
 (X)_{I*J}(a_0)\ar[u]^\simeq\ar[rr]^-\simeq &  
 & 
(X)_{I*J}(b_1)\ar[u]^\simeq\\
X_I(a_0)\ar[u]^{\simeq}\ar[r]^-{\simeq} &
X_I(a_1)\cong (X)_{I*J}(a_1)= (X)_{I*J}(b_0)\cong X_J(b_0)\ar[r]^-{\simeq}&
X_J(b_1),\ar[u]^\simeq
}\]
where we use the notation $X_Z=X\times_Y Z$, the vertical maps in the upper square are induced by base change of $w*v$,
and the isomorphisms $X_I(a_i)\cong X_{I*J}(a_i)$ and $X_J(b_i)\cong X_{I*J}(b_i)$, for $i=0,1$,
are induced by the natural closed immersions $I\inj  I*J$ and $J\inj I*J$.
The upper square commutes by the definition of $\mu_{w*v}$, the lower square commutes by the construction 
of the isomorphisms, see the proof of Proposition \ref{prop:monodromy}.
Thus the whole square commutes and by definition of the monodromy we obtain
the equality 
\eq{para:loop-group1}{\mu_{w*v}=\mu_{v}\circ \mu_w}
of automorphisms of the fiber functor  $\Phi_{y_0}:(\FinEt/Y)\to (\FinSet)$.
Furthermore if $u: (I,a_0, a_1)\to (Y, y_0, y_0)$, $v: (J,b_0, b_1)\to (Y, y_0, y_0)$, and 
 $w: (L,c_0, c_1)\to (Y, y_0, y_0)$ are three loops based at $y_0$, then  the universal property of pushout diagrams yields 
 a canonical isomorphism $\tau: I*(J*L)\to (I*J)*L$ and  an equality
\eq{para:loop-group2}{(u*v)*w\circ \tau = u*(v*w)\quad  \text{in } (\Sch/Y).}

Clearly, homotopy  between paths defines an equivalence relation $w\sim v$  and we denote by 
$$\pi_0\Omega^\alg(Y,y_0)$$ 
the set of homotopy classes of algebraic loops based at $y_0$.
We denote by $[w]$ the homotopy class of a loop $w:I\ra Y$ at $y_0$.  
According to \eqref{para:loop-group1} we obtain a well defined operation
\[*: \pi_0\Omega^\alg(Y,y_0)\times \pi_0\Omega^\alg(Y,y_0)\lra \pi_0\Omega^\alg(Y,y_0),\quad
([w], [v])\longmapsto [w*v].\]
This operation is associative by $\eqref{para:loop-group2}$ and clearly 
any constant loop $I\to y_0\to Y$ has the same  homotopy class denoted by $e$, which defines a neutral element for
$*$. Moreover,  for a loop $w: (I, a_0, a_1)\to (Y, y_0, y_0)$  we define the loop
$w': (I, a_1, a_0)\to (Y,y_0,y_0)$ by switching the end points of $I$. Clearly
$[w]*[w']=e$. Hence concatenation of algebraic loops defines a group structure on $\pi_0\Omega^\alg(Y,y_0)$,
which we therefore call the {\em algebraic loop group}.

By definition of the algebraic fundamental group, our homotopy relation, and the relation
\eqref{para:loop-group1} the algebraic loop space $\Omega^\alg(Y,y_0)$
induces an injective homomorphism
\eq{homotopy classes1}{
\pi_0\Omega^\alg(Y,y_0)^{\op}\lra \pi_1^\alg(Y,y_0),\quad w\longmapsto \mu_w
}
of groups, where we use the opposite group structure on the left hand side.
We  regard this as an inclusion of groups.
\end{para}

The following is the main result of this note.

\begin{theorem}\label{thm:loop-et}
Let $X$ be a  connected scheme that is separated and of finite type over a field $k$, endowed  with a geometric point  $x_0:\Spec(k^\sep)\to X$.
Then the canonical injection \eqref{homotopy classes1} has dense image.
It is actually bijective, provided   $X$ is proper.
\end{theorem}

The proof of Theorem \ref{thm:loop-et} requires some preparations and will be given in Section \ref{sec:proof}.
We remark  that we do not expect \eqref{homotopy classes1} to be an isomorphism for non-proper schemes.


\section{A Lefschetz type theorem}
\mylabel{General Lefschetz}

The   \emph{Lefschetz Hyperplane Theorem} gives a strong relation between  the homology of 
a  projective   complex manifold $X$ of dimension $n\geq 2$
and the homology of an   ample divisor $D\subset X$. 
The original arguments appear in \cite{Lefschetz 1924}, Chapter V, Section III.
Analogous statements for  fundamental groups were first obtained 
by Bott \cite{Bott 1959}: The induced map 
$\pi_1(D,x_0)\ra \pi_1(X,x_0)$   is bijective provided $n\geq 3$, and at least surjective if $n\geq 2$.

The latter statement extends  to projective schemes $X$ over arbitrary ground fields $k$:
According to  \cite{SGA 2}, Expos\'e XII, Corollary 3.5 the map $\pi^\alg_1(D,x_0)\ra \pi^\alg_1(X,x_0)$ is surjective
provided that $\depth(\O_{X,a})\geq 2$ for each closed point $a\in X$.
If $X$ is additionally   \emph{Cohen--Macaulay}, 
i.e., at every point $a\in X$ the equality $\dim(\O_{X,a})=\depth(\O_{X,a})$ holds, the above can be iterated
and one finds  a connected curve $C\subset X$ such that  $\pi^\alg_1(C, x_0)\to \pi^\alg_1(X, x_0)$ is surjective.

In this section we generalize the latter statement to arbitrary proper $k$-schemes.

\begin{para}\mylabel{curve property}
Let $X$ be a non-empty connected noetherian scheme. 
We consider the following property:

\begin{enumerate}[label={(C)}]
\item\label{C}   For each closed subscheme $Z\subset X$ with $\dim(Z)\leq 0$, there is a  connected  closed  subscheme
$C\subset X$  with $0\le \dim(C)\le 1$ and $Z\subset C$ such that, for each finite \'etale covering $U\ra X$ with connected total space, the restriction $C_U=C\times_X U$ remains connected. 
\end{enumerate}
\end{para}

\begin{remark}\label{rmk:C}
We remark that  property \ref{C} does not hold for   affine schemes in general, as 
the following simple example  shows (confer Lemma 5.4 in \cite{Esnault 2017}):
Let $k$ be a ground field of characteristic $p>0$, and $C$ be an connected affine  plane curve inside
$\AA^2=\Spec(R)$, defined by some non-constant polynomial $f=f(x,y)$ inside the ring $R=k[x,y]$.
Then  there exists a connected finite \'etale covering  $X\to \AA^2$
whose restriction to $C$ becomes trivial. To see this, take any  $h\in (f)$   not of the form $g^p-g$ with $g\in R$.
Via the identification 
$$
H^1_{\et}(\AA^2_k, \ZZ/p)=R/\{ g^p-g\mid g\in R\}
$$
coming from the  Artin--Schreier sequence,  the polynomial
$h$ corresponds to a non-trivial $\ZZ/p$-torsor $X\to \AA^2_k$,
in particular it is a connected finite \'etale covering of $\AA^2_k$.
On the other hand $h$ maps to zero in $\bar{R}/\{a^p-a\mid a\in \bar{R}\}$, where $\bar{R}=R/(f)$. In other words, the restriction of 
$X$ to $C$ is trivial.
\end{remark}

\begin{lemma}
\mylabel{reduction step}
Let $X$ be a non-empty connected noetherian scheme,   $f:X'\ra X$ be a proper and surjective  morphism,
and $X'_v\subset X'$ be the connected components.  
If property  \ref{C} holds for all  $X'_v$, then it also holds for $X$.
\end{lemma}

The proof of this lemma is inspired by the Seifert--van Kampen Theorem from \cite{Stix 2006}, but  is more elementary.
Note that  even if one wants to use  property \ref{C} on $X$ with $Z=\varnothing$, the proof for the lemma relies in an essential way on property \ref{C} 
on $X'_v$  with non-empty $Z'$. We first gather some  basic material on graphs.

\begin{para}\label{para:graph}
Let $f:X'\to X$ be as in the statement of Lemma \ref{reduction step}.
Set  $X''=X'\times_XX'$ and write $\pr_1,\pr_2:X''\ra X'$  for the two projections.
The schemes $X'$ and $X''$  are noetherian,  hence the sets of connected components $\pi_0(X')$ 
and $\pi_0(X'')$ are finite.
Consider the induced maps
$$
\pr_1\times\pr_2:\pi_0(X'')\lra \pi_0(X')\times\pi_0(X').
$$
This defines an \emph{oriented  graph} $\Gamma=(E,V,\pr_1\times\pr_2)$
 in the sense of Serre \cite{Serre 1980}, Section 2.1:
the set of \emph{vertices} is  $V=\pi_0(X')$, and the set of \emph{oriented edges} is $E=\pi_0(X'')$.
The endpoints of an edge $e\in E$ are the images $v_i=\pr_i(e)$. The orientation is given by
declaring $v_1$ as the initial vertex, and $v_2$ as the terminal vertex. 
We usually write $X'_v\subset X'$ and $X''_e\subset X''$ for the connected components corresponding to a vertex $v$
and an edge $e$.

Note that edges could have the same initial and terminal vertices, and 
several edges could share their initial and  terminal vertices.
By abuse of notation, we also write $v\in \Gamma$ and $e\in\Gamma$ to 
denote vertices and edges of the graph, if there is no risk of confusion.
A morphism $f:\Gamma\ra \Gamma'$ between oriented graphs comprises compatible maps $V\ra V'$ 
and $E\ra E'$. We simply say that $f$ is a \emph{map of oriented graphs}.
Also note that the graph $\Gamma$ constructed above is connected, since the scheme $X$ is connected.
\end{para}

\begin{proof}[Proof of Lemma \ref{reduction step}.]
Since $f:X'\to X$ is proper,  the image $f(Z)$ of a closed subscheme $Z\subset X'$ is closed and satisfies  $\dim f(Z)\le \dim (Z)$.
In particular, closed points are mapped to closed points. 
The assertion is trivial for $\dim(X)\le 1$. We now  assume $\dim X\ge 2$. We use the notation from \ref{para:graph}.
Let $Z\subset X$ be a zero-dimensional closed subscheme or the empty set.
For each edge $e\in \Gamma$ choose a closed point $x_e\in X''_e$.
Set
\[x_{e,i}:= \pr_i(x_{e})\in X'_{v_i},\quad \text{where } v_i=\pr_i(e)\in \Gamma, \, i=1,2.\]
Since $\Gamma$ is finite we find for each vertex $v\in \Gamma$ a $0$-dimensional closed subset $Z'_v\subset X'_v$,
such that 
\begin{enumerate}[label=(\alph*)]
\item\label{red-step-a} $Z\cap f(X_v')\subset f(Z'_v)$ and 
\item\label{red-step-b}  $x_{e,i}\in Z'_v$, for all edges $e\in\Gamma$ with  $v=\pr_i(e)$ for $i=1$ or $2$.
\end{enumerate}
Condition \ref{red-step-b} is immediate, and one can achieve  \ref{red-step-a} by picking a closed point in each of the finitely many schemes $f^{-1}(z)\cap X'_v$,
with  $z\in Z\cap f(X_v')$.
By the surjectivity of $f$ we have
\eq{red-step1}{Z\subset \bigcup_v f(Z'_v).}
Applying \ref{C} to $X'_v$ and $Z'_v$ we find an 
at most 1-dimensional closed subscheme  $C'_v\subset X'_v$ containing $Z'_v$,
such that the pullback of any connected finite \'etale covering of $X'_v$ to $C_v'$ stays connected. 
Then $C=\bigcup_v f(C'_v)$ is closed, at most 1-dimensional, and contains $Z$, by \eqref{red-step1}.
It remains to show that for each  connected finite \'etale covering $U\to X$ the pullback $U\times_X C$  remains connected
(then $C$ is connected as well).

To this end fix such a finite \'etale covering $u:U\to X$, with $U$ non-empty and connected.
Denote by $\Gamma_U$ the graph defined by $\pr_1\times\pr_2:\pi_0(U'')\to \pi_0(U')\times\pi_0(U')$, 
where $U'=U\times_X X'$ and $U''=U\times_{X}X''=U'\times_U U'$.
We obtain a surjection of graphs $u:\Gamma_U\to\Gamma$ 
and for each edge $\epsilon\in \Gamma_U$ we obtain a finite and \'etale morphism between connected schemes 
$U''_{\epsilon}\to X''_{u(\epsilon)}$ which therefore is surjective.
Thus for each edge $e\in \Gamma$ and edge $\epsilon\in\Gamma_U$  mapping to $e$ we can choose a  
closed point $x_{U,\epsilon}\in U_{\epsilon}''$ with $u(x_{U,\epsilon})=x_e$.

For a vertex $w\in \Gamma_U$ mapping to $v\in \Gamma$ denote by $I_{w}$  the image of $U'_{w}\times_{X'_v} C'_{v}$ under the map
\[U'\times_{X'} C'_{v}=U\times_X C'_{v}\lra U\times_X C_{v},\quad \text{where } C_v=f(C'_v).\]
The map is induced by the base change with the composition 
$C'_v\inj X'_v\inj X'\xr{f} X$,  and therefore is closed and surjective. 
Hence 
\eq{red-step2}{U\times_X C_{v}=\bigcup_{w\in u^{-1}(v)} I_{w} \quad \text{and}\quad 
U\times_X C=\bigcup_{w \text{ edge in } \Gamma_U} I_w,}
where each $I_{w}$ is closed. By our choice of $C'_v$ the pullback of the connected \'etale covering 
$U'_w\to X'_v$ over $C'_v$ remains connected. Thus $I_w$ is the image of a connected scheme and is hence connected.
Let $w_1$ and $w_2$ be  the initial and the terminal vertices of an edge $\epsilon\in \Gamma_U$, then 
$x_{U,\epsilon}\in U''_\epsilon$ maps via the $i$th projection to points
$\pr_i(x_{U,\epsilon})$ in $U'_{w_i}\times_{X'_{u(w_i)}} C'_{u(w_i)}$, $i=1,2$, and these points map  to same point 
in  $U$.
Thus the intersection $I_{w_1}\cap I_{w_2}$ is non-empty, if $w_1$ and $w_2$ are linked by an edge in $\Gamma_U$. 
Since the graph $\Gamma_U$ is connected so is the scheme  $U\times_X C$.
This completes the proof.
\end{proof}

\begin{proposition}\label{prop:C}
\mylabel{C over a field}
Let $X$ be a connected scheme that is proper over a field $k$.
Then  $X$ has property $\ref{C}$.
\end{proposition}

\begin{proof}
We  proceed by induction on  $n=\dim X$.
There is nothing to prove for $n\le 1$.
Assume $n\ge 2$ and that \ref{C} holds for all  connected schemes that are proper over $k$ 
and have dimension  $\le n-1$. 
Using Lemma \ref{reduction step}   
we can make the following reductions:
\begin{enumerate}
\item $X$ reduced (using the proper bijection $X_\red\to X$);
\item $X$ projective over $k$ (using Chow's Lemma);
\item $X$ integral (using the proper surjection $\coprod X_i\to X$, with $X_i$ the irreducible
 components of $X$).
\end{enumerate}
According to Kawasaki's result (\cite{Kawasaki 2000}, Theorem 1.1) there is a proper birational $X'\ra X$ 
such that the scheme $X'$ is Cohen--Macaulay. Moreover, this Macaulayfication   arises as a sequence of blowing-ups.
Hence applying Lemma \ref{reduction step} one more time, we are reduced to consider the case
that $X$ is projective, integral, and Cohen--Macaulay over a field $k$.  
Let $Z\subset X$ be a closed subscheme with $\dim(Z)\le 0$.
By, e.g., \cite{Gabber Liu Lorenzini 2015}, Theorem 5.1, we find  an effective ample divisor  $D\subset  X$ containing $Z$.
By induction the following claim implies that $X$ satisfies \ref{C}:
\begin{claim}
Let  $X'\ra X$ be an \'etale covering whose total space is connected.  
Then the restriction $D'=X'\times_X D$ remains connected. 
\end{claim}
The   argument to prove the claim is similar to \cite{Hartshorne 1970}, Chapter II, Corollary 6.2.
Let us recall it for the sake of completeness:
Consider the ample invertible sheaf  $\shL=\O_X(D)$.
Since $X'\ra X$ is finite and surjective, the inclusion $D'\subset X'$ remains an effective Cartier divisor,
and the corresponding invertible sheaf is the pullback $\shL'=\shL|X'$, which is still ample. 
Since $X$ is projective and Cohen--Macaulay, so is $X'$. 
Let $\omega_{X'}$ be the dualizing sheaf over $k$.
Then $h^1(\shL'^{\otimes -t}) = h^{n-1}(\omega_{X'}\otimes\shL'^{\otimes t})$
for every integer $t$. The right hand side vanishes for $t$ sufficiently large, because $\shL'$ is ample and $n-1\geq 1$.
Replacing $D$ by $tD$, we may assume $H^1(X',\shL'^{\otimes-1})=0$.
The short exact sequence 
$0\ra \shL^{\otimes-1}\ra\O_{X'}\ra\O_{D'}\ra 0$ thus gives a surjection of rings
$
H^0(X',\O_{X'})\ra H^0(D',\O_{D'}).
$
The term on the left is a finite extension of the ground field $k$ because $X'$ is integral and proper.
Hence the above map is bijective, and $D'$ must be connected.
\end{proof}

In view of \ref{para:functorGalois} we obtain the following corollary.
\begin{corollary}\label{cor:C}
Let $k$ be a field and set $K=k^\sep$.
Let $X$ be a connected scheme which is proper over $k$ and let 
$x_0:\Spec (K)\to X$ be a geometric point.
Then there exists a connected, reduced, affine, and 1-dimensional scheme $C$ of finite type over $k$ 
and a $k$-morphism $C\to X$, such that $x_0$ factors via $C$ and the induced map 
\[\pi_1^\alg(C, x_0)\lra \pi_1^\alg(X, x_0)\]
is surjective.
\end{corollary}
\begin{proof}
By Proposition \ref{prop:C} and \ref{para:functorGalois} we  find a connected  closed subscheme 
$C_1\subset X$ of dimension at most 1, such that $x_0$ factors via $C_1$ and the induced map 
\[\pi_1^\alg(C_1, x_0)\lra \pi_1^\alg(X, x_0)\]
is surjective. Since passing to the reduced subscheme does not change the fundamental group, 
we may assume  $C_1$ reduced. If $\dim (C_1)=0$, then $C_1=\Spec(L)$ with $L$ a subfield of  $K$.
In this case we can take $C:=\AA^1_L$ with map $\AA^1_L\to \Spec (L)=C_1\to X$ and 
factorization of $x_0$ given by the composition of $\Spec (K)\to \Spec (L)$ with the inclusion of the zero section
$\Spec(L)\inj \AA^1_L$.

Assume $\dim (C_1)=1$.  Note that $C_1$ is quasi-projective, hence we find an affine open $U\subset C_1$
which is connected and contains the singular locus of $C_1$ and the image point of $x_0$. 
In particular we remove from $C_1$ only finitely many regular closed points, whose local rings 
are therefore discrete valuation rings. 
Hence it follows from \cite[\href{https://stacks.math.columbia.edu/tag/0BSC}{Tag 0BSC}]{stacks-project}
that $\pi_1^\alg(U,x_0)\to \pi_1^\alg(C_1,x_0)$ is surjective and we can take $C=U$.
\end{proof}


\section{The non-proper case}

In view of Remark \ref{rmk:C}, we consider the following weaker variant of condition \ref{C} in this  section.
\begin{para}\label{para:C*}
Let $X$ be a non-empty connected noetherian scheme. 
We consider the following  property:
\begin{enumerate}[label=($\text{C}^*$)]
\item\label{C*} For each closed subscheme $Z\subset X$ with $\dim(Z)\leq 0$ and each finite 
\'etale covering $U\to X$ with connected total space there exists  a  connected   closed  subscheme
$C\subset X$  with $0\le \dim(C)\le 1$ and $Z\subset C$, such that $C_U=C\times_X U$ remains connected.
\end{enumerate}
\end{para}

\begin{lemma}
\mylabel{reduction step*}
Let $X$ be a non-empty connected noetherian scheme. 
Let $f:X'\ra X$ be a  universally closed and surjective  morphism from a noetherian scheme $X'$, with connected components $X'_v$.
If property \ref{C*} holds for all $X'_v$, then it also holds for $X$.
\end{lemma}
We remark, that  besides proper maps, all integral morphisms  are universally closed, see \cite[Proposition (6.1.10)]{EGAII}.
Thus the above lemma applies to the normalization map  $X'\ra X$ of an integral scheme,
and also to the projection from the base-change $X'=X\otimes_kk'$ with respect to any algebraic ground field extension, provided that $X'$ stays noetherian.

Since in Lemma \ref{reduction step*} the map $f:X'\to X$ is not assumed to be of finite type,
the product $X'\times_X X'$ may not be noetherian and might have infinitely many connected components.
Thus the graph constructed in \ref{para:graph} might have an infinite set of edges.
To deal with this we record the following lemma.

\begin{lemma}\label{lem:finite-graph}
Let $\Gamma$ be a connected  oriented graph with a finite set of vertices.
Then there exists a finite oriented subgraph $\Gamma'\subset \Gamma$, which has the same set of vertices as $\Gamma$ and which is connected.
\end{lemma}
\begin{proof}
Choose for each pair of vertices  $v,w\in \Gamma$ a path $p_{v,w}$ connecting them. Take $\Gamma'$ to be the graph whose set of vertices is equal to the set of vertices of  $\Gamma$ and
whose edges are given by all the finitely many edges appearing in the paths $p_{v,w}$, for all pairs of vertices $(v,w)$. 
\end{proof}

\begin{proof}[Proof Lemma \ref{reduction step*}.]
Since $f$ is closed and surjective  the image  $f(Z)$ of a closed subscheme $Z\subset X'$ 
is closed and satisfies $\dim f(Z)\le \dim (Z)$, see \cite[Proposition (5.4.1)]{EGAIV2}.
In particular, closed points are mapped to closed points.

We assume $\dim X\ge 2$.
Let $Z\subset X$ be a closed subset with $\dim(Z)\le 0$. 
Set $X''=X'\times_X X'$. Let $u:U\to X$ be a finite \'etale covering with connected total space.
We want to find a closed connected at most one dimensional subscheme $C$ which contains $Z$
and such that the pullback of $U$ over $C$ stays connected. To this end we may assume that 
$u: U\to X$ is a finite \'etale Galois covering with Galois group $G$.
We denote by $U'=U\times_X X'$ and $U''= U\times_X X''$ the base changes and by 
$\Gamma_X$ and $\Gamma_U$ the oriented graphs defined as in  \ref{para:graph} by 
${\pr_1}\times\pr_2:\pi_0(X'')\to \pi_0(X')\times \pi_0(X')$ and  ${\pr_1}\times\pr_2:\pi_0(U'')\to \pi_0(U')\times \pi_0(U')$,
respectively. These graphs are connected, since $X$ and $U$ are, and have finite sets of vertices. 
The map $u$ induces a surjective map of graphs $\Gamma_U\to \Gamma_X$ again denoted by $u$.
For any vertex $w\in\Gamma_U$ with $v=u(w)\in \Gamma_X$, the morphism $u$ induces a 
finite \'etale morphism $U'_w\to X'_v$.

Let  $\Gamma'_U\subset \Gamma_U$ be a finite connected subgraph with the same vertices as $\Gamma_U$, 
see Lemma \ref{lem:finite-graph}. 
Choose closed points $x_{U,e}\in U''_e$ for any edge  $e\in \Gamma'_U$.  
Set
\[x_{e,i}:= u(\pr_i(x_{U,e}))\in X'_{v_i},\quad \text{where } v_i=u(\pr_i(e))\in \Gamma, \, i=1,2.\]
Since $\Gamma'_U$ is finite we find as in the proof of Lemma \ref{reduction step} 
for each vertex $v\in \Gamma_X$ a $0$-dimensional closed subset $Z'_v\subset X'_v$ such that 
\begin{enumerate}[label=(\alph*)]
\item\label{red-step*-a} $Z\cap f(X_v')\subset f(Z'_v)$ and 
\item\label{red-step*-b}  $x_{e,i}\in Z'_v$, for all edges $e\in\Gamma'_U$ with  $v=u(\pr_i(e))$ for $i=1$ or $2$.
\end{enumerate}
By the surjectivity of $f$ we have
\eq{red-step*1}{Z\subset \bigcup_v f(Z'_v).}
Fix a vertex $v\in \Gamma$ and choose $w_0\in \Gamma'_U$ mapping to $v$. 
Applying \ref{C*} to the finite \'etale covering $U'_{w_0}\to X'_v$ we find an at most 1-dimensional 
connected closed subscheme $C'_{v}\subset X'_v$ containing $Z'_v$, 
such that the restriction $U'_{w_0}\times_{X'_v} C'_{v}$ remains connected.
The base change $U\times_X C'_v\to C'_v$ is a Galois covering with Galois group  $G$.
In particular $G$ acts transitively on the connected components of $U\times_X C'_v$ and 
we obtain isomorphisms 
\eq{red-step*2}{U'_{w}\times_{X'_v} C'_v\cong U'_{w_0}\times_{X'_v} C'_v,\quad 
\text{for all }w\in \Gamma_U \text{ mapping to } v,}
in particular all these schemes are connected.
Set
\[C_{v}:=f(C'_{v})\quad \text{and}\quad C:=\bigcup_{v} C_{v}.\]
It follows that $C\subset X$ is closed, non-empty and at most 1-dimensional.
By  \eqref{red-step*1} and  $Z'_v\subset C'_{v}$ we have $Z\subset C$.
Moreover, using the choice of $C'_v$ together with \eqref{red-step*2} 
we can argue in the same way as in the last paragraph of the proof of Lemma \ref{reduction step} with $\Gamma_U$ there replaced by
$\Gamma'_U$ here to deduce that  $U\times_X C$   is connected. 
\end{proof}

\begin{lemma}\mylabel{Bertini}
Let $k$ be an algebraically closed field and $X$ an integral  quasi-projective $k$-scheme of $\dim X\ge 2$.
Let $Z\subset X$ be a finite (possibly empty) set of closed points and 
$X'\to X$  a finite and surjective morphism with $X'$ irreducible. 
Then there is an integral  closed subscheme  $H\subset X$ of codimension $1$ containing $Z$, 
such that $X'\times_X H$ is irreducible.
\end{lemma}

\begin{proof}
This follows from a classical Bertini theorem, where we use a trick of Mumford to ensure that the hyperplane
contains $Z$, see the proof of  the Lemma on p.\ 56 in \cite{Mumford 1970}:
Let  $f: Y\to X$ be the blowing up with center $Z$.
Then   $\dim f^{-1}(z)\ge 1$, for all $z\in Z$. 
We fix an embedding $Y\inj \PP^n_k$. Denote by $Y'\subset Y\times_X X'$ an irreducible component
which maps birationally onto $X'$.
By \cite[I, Corollary 6.11, 3)]{Jouanolou 1983} applied to the quasi-finite  morphism 
$Y'\to Y\to \PP^n_k$ we find a hyperplane $H_1\subset \PP^n_k$ 
which is not contained in the exceptional locus
of $f$ and such that its pullback to $Y'$ is irreducible. (Here we use $k$ algebraically closed, since in {\em loc. cit.} $Y'$ is required to be geometrically irreducible and 
$k$ to be infinite.) Since $f^{-1}(z)$ is closed in $\PP^n_k$ and 
$\dim f^{-1}(z)+\dim H_1\ge n$ we find $f^{-1}(z)\cap H_1\neq \emptyset$, for all $z\in Z$.
Set $H:=f(H_1\cap Y)_\red$. Then $Z\subset H$ 
and the pullback of $H$ to $X'$ is birationally dominated by $H_1\times_{\PP^n} Y'$ and hence is irreducible.
\end{proof}

\begin{proposition}\label{C* over a field}
Let $X$ be a connected scheme which is separated and of finite type  over a field $k$.
Then  $X$ has property \ref{C*}.
\end{proposition}

\begin{proof}
We proceed by induction on  $d=\dim X$. There is nothing to prove for $d\le 1$.
Assume $d\ge 2$ and that \ref{C*} holds for all  connected schemes that are separated and of finite type over $k$ 
and have dimension  $\le d-1$. 
Since $X$ is connected all its irreducible components have dimension $\ge 1$.
Using Lemma \ref{reduction step*}   
we can therefore make the following reductions:
\begin{enumerate}
\item $k$ algebraically closed (since for $\bar{k}$ the algebraic closure of $k$ the morphism $X\otimes_k\bar{k}\to X$ is integral, whence universally closed);
\item $X$ reduced (by considering the proper morphism $X_\red\to X$);
\item $X$ quasi-projective over $k$ (Chow's Lemma);
\item $X$ integral (by considering the proper and surjective morphism $\sqcup X_i\to X$, with $X_i$ the irreducible components of $X$);
\item $X$ normal (by considering the normalization $\widetilde{X}\to X$).
\end{enumerate}
Assume $k$ and $X$ satisfy the conditions above. 
Let $U\to X$ be a finite \'etale map with $U$ connected. By the normality of $X$ the scheme $U$ is normal as well.
Thus $U$ is irreducible. Hence the existence of a curve for $U\to X$ as in \ref{C*} follows directly from Lemma \ref{Bertini}
and the induction hypothesis.
\end{proof}

\begin{corollary}\label{cor:C*}
Let $k$ be a field and set $K=k^\sep$.
Let $X$ be a connected scheme which is separated and of finite type over $k$ and let 
$x_0:\Spec (K)\to X$ be a geometric point.
Let $X'\to X$ be a finite \'etale Galois covering.
Then there exists a connected, reduced, affine, and 1-dimensional scheme $C$ of finite type over $k$ 
and a $k$-morphism $C\to X$, such that $x_0$ factors via $C$ and the composite map 
\eq{cor:C*1}{\pi_1^\alg(C, x_0)\lra \pi_1^\alg(X, x_0)\lra\Aut(X'/X)^{\rm op},}
is surjective. Here the second map is the natural surjection from \cite[Exp. V, 4, h)]{SGA 1}.
\end{corollary}
\begin{proof}
In general the composition \eqref{cor:C*1} is surjective if the pullback of $X'$ over $C$ stays connected.
Hence the statement follows from Proposition \ref{C* over a field} the same way Corollary \ref{cor:C} follows from
Proposition \ref{prop:C}.
\end{proof}

\section{Proof of the main theorem}\label{sec:proof}

We prove Theorem \ref{thm:loop-et}. We  start by proving the second statement, i.e.,
for  a  {\em proper} and connected scheme $X$ over a field $k$ with geometric point $x_0:\Spec (k)\to X$ we want to show that
the natural injective group homomorphism $\pi_0\Omega^\alg(X,x_0)\to \pi_1^\alg(X, x_0)$ is surjective as well.

Set $K:=k^\sep$. By Corollary \ref{cor:C}
we find a connected, affine, reduced and 1-dimensional $k$-scheme of finite type $C$ with a  morphism
$C\to X$ such that $x_0$ factors via $C$ and the natural 
$\pi_1^\alg(C, x_0)\surj \pi_1^\alg(X,x_0)$ is surjective.
We obtain a commutative diagram
\eq{sec:proof1}{\xymatrix{
\pi_0\Omega^\alg(C,x_0)^{\op}\ar[r]^{}\ar[d] & \pi_1^\alg(C,x_0)\ar@{->>}[d]\\
\pi_0\Omega^\alg(X,x_0)^{\op}\ar[r]^{} & \pi_1^\alg(X,x_0),
}}
in which the vertical arrow on the right is surjective.
It therefore remains to show that the top horizontal arrow is surjective.
To this end we choose some pro-object $(C_i)_{i\in I}$ in $(\FinEt/C)$ representing $\Phi_{x_0}$, see \cite{SGA 1}, Expos{\'e} V, 4. 
We write ${\rm Pro}(\FinEt/C)$ for pro-objects in $(\FinEt/C)$.
Since $I$ is a filtered set and the transition maps  $C_i\to C_j$ are finite \'etale, and therefore affine,
we may form the projective limit $\widetilde{C}=\varprojlim_{i\in I} C_i$ in the category of $C$-schemes,
see \cite{EGAIV3}, Proposition (8.2.3).
For a $C$-scheme $T$ we obtain a functorial isomorphism
\[\Hom_C(T, \widetilde{C})=\varprojlim_i \Hom_C(T, C_i).\]
In particular an element $a_0\in \varprojlim \Phi_{x_0}(C_i)$
is a morphism $a_0: \Spec (K)\to \widetilde{C}$ over $x_0$.
Furthermore, by \cite{SGA 1}, Expos{\'e} V, 4, h), we have
\begin{align*}
\Aut_{{\rm Pro}(\FinEt/C)}((C_i)_{i}) &=\Hom_{{\rm Pro}(\FinEt/C)}((C_i)_i, (C_j)_j) \\
                                                         &= \varprojlim_j \Hom_{{\rm Pro}(\FinEt/C)}((C_i)_i, C_j) \\
                                                         &=\varprojlim_j \Hom_C(C_j, C_j)
                                                         = \varprojlim_j \Hom_C(\widetilde{C}, C_j) \\
                                                         &= \Hom_C(\widetilde{C}, \widetilde{C}) = \Aut_C(\widetilde{C}).
\end{align*}
By {\em loc. cit.}, the choice of $a_0: \Spec (K)\to \widetilde{C}$ over $x_0$ yields a functorial isomorphism
\[\Hom_C(\widetilde{C}, U)\stackrel{\simeq}{\lra} \Phi_{x_0}(U) , \quad f\longmapsto f\circ a_0,\qquad U\in (\FinEt/C).\]
We obtain an  isomorphism
\[\theta:\Aut_C(\widetilde{C})^{\rm op}\stackrel{\simeq}{\lra} \Aut(\Phi_{x_0})=\pi_1^\alg(C, x_0),\]
which sends a $C$-automorphism $\sigma: \widetilde{C}\to \widetilde{C}$
to the automorphism $\theta(\sigma)$ of $\Phi_{x_0}$, which on $U\in (\FinEt/C)$ is given by
\[\Phi_{x_0}(U)\ni f\circ a_0\longmapsto f\circ (\sigma\circ a_0)\in \Phi_{x_0}(U), \quad f\in \Hom_C(\widetilde{C}, U).\]
We claim that for any $\sigma$ the automorphism $\theta(\sigma)$ is in the image of 
$\pi_0\Omega^\alg(C, x_0)^{\op}$.
Indeed, by construction  $\widetilde{C}$ is affine, reduced, connected and 1-dimensional and is simply connected
We have the $K$-rational point $a_0: \Spec (K)\to \widetilde{C}$ over 
$x_0: \Spec (K)\to C$. We note that $a_0$ is a closed immersion. Indeed, for any  finite separable field extension $L/k$
a connected component $C_0$ of $C\otimes_k L$ is a finite \'etale covering of $C$; hence we have a map $\widetilde{C}\to C_0$.
It follows that the algebraic closure of $k$ in  $H^0(\widetilde{C}, \O_{\widetilde{C}})$ is equal to $K=k^{\sep}$, which implies that 
$a_0$ is a closed immersion. 
Thus $(\widetilde{C},a_0, \sigma\circ a_0)$ is an interval scheme in the sense of
Definition \ref{interval scheme} and the map $\widetilde{C}\to C$ induces an algebraic loop
$w: (\widetilde{C}, a_0, \sigma\circ a_0)\to (C, x_0, x_0)$. The map \eqref{homotopy classes1}
sends the loop $w$ to the monodromy $\mu_w$ which by construction is equal to $\theta(\sigma)$.
This completes the proof of the second part of the theorem.

It remains to show that assuming  $X$ is connected and only separated and of finite type over $k$, then 
the image of $\pi_0\Omega^\alg(X,x_0)^{\op}\to\pi_1^\alg(X,x_0)$ is dense. It suffices to show that for any finite \'etale Galois covering
$X'\to X$ the composition
\[\pi_0\Omega^\alg(X,x_0)^{\op}\lra\pi_1^\alg(X,x_0)\lra \Aut(X'/X)^{\rm op}\]
is surjective. By Corollary \ref{cor:C*} we find a curve $C$ as above such that 
the composition 
\[\pi_1^\alg(C,x_0)\lra \pi_1^\alg(X,x_0)\lra \Aut(X'/X)^{\rm op}\]
is surjective. Thus the statement follows from the surjectivity of top horizontal map in \eqref{sec:proof1}
proved above.
\qed

\end{document}